# Mixed finite element formulation in large deformation frictional contact problem

**Laurent Baillet\* — Taoufik Sassi\*\***

*\* Laboratoire de Mécanique des Contacts et des Solides, CNRS-UMR 5514*
*INSA de Lyon, 69621 Villeurbanne Cedex, France*
*e-mail: laurent.baillet@insa-lyon.fr*

*\*\* Laboratoire de Mathématiques Nicolas Oresme*
*CNRS 6139 Université de Caen*
*Bd du Marechal Juin, 14032 CAEN Cedex*
*e-mail sassi@math.unicaen.fr*

ABSTRACT. *This paper presents a mixed variational framework and numerical examples to treat a bidimensional friction contact problem in large deformation. Two different contact algorithms with friction are developed using the 2D finite element code PLAST2[1].*
*The first contact algorithm is the classical node-on-segment, and the second one corresponds to an extension of the mortar element method[2,3,4] to a unilateral contact problem with friction. In this last method, the discretized normal and tangential stresses on the contact surface are expressed by using either continuous piecewise linear or piecewise constant Lagrange multipliers in the saddle-point formulation. The two algorithms based on Lagrange multipliers method are developed and compared for linear and quadratic elements.*

KEY WORDS: contact, mixed finite element, friction, dynamic explicit, mortar elements.

## 1. Introduction

In order to solve a contact problem with friction, it is necessary to possess numerical tools adapted to the strong non-linearity, the incompatibility of meshing on the contact zone and the evolutionary characteristics of the surfaces.

Methods used differ by their contact algorithm, their time integration scheme and the construction of the global contact matrices obtained from the writing of the non interpenetration condition. The solution to a contact problem is obtained using various methods such as the penalisation or the Lagrange multiplier methods[5,6,7,8,9,10]. Among these last methods one finds the gradient methods[11,12], those of the increased Lagrangian or other mixed approaches[13,14,15].

For the resolution of a contact problem without friction using the Lagrange multipliers method, the construction of the global contact matrices used for the calculation of contact stresses is not unique. Usually, the classical node-on-segment strategy (local type approach) is used. The mortar element method[3] initially presented for domain decomposition has been used for the resolution of a unilateral contact problem[16,17,18]. The introduction of another Lagrange multiplier related to the tangential friction stress following a Coulomb or Tresca law has been presented by the researchers[4,19]. McDevitt and Laursen have used the mortar element method in the case of small deformations, by using the penalisation method and for problems of non evolutionary contact surfaces.

In this paper, the contact treatment is based on the mortar-finite elements method (global type approach), it uses the Lagrange multiplier method and is developed for large deformation problems where the contact surface is evolutionary. These algorithms have been implemented and tested in the 2D finite element code PLAST2. PLAST2 includes large deformation and non linear material behavior. It is based on a Lagrangian-mesh Cauchy-stress formulation in conjunction with an explicit time integration scheme. The forward Lagrange multiplier method is used to treat the contact between deformable bodies.

The paper is organized as follows. First we introduce the equations modeling the Signorini problem with friction, continuous mixed variational formulation and contact algorithms are presented. Then several numerical simulations which include contact between deformable bodies are performed. Comparisons of the two contact algorithms, the choice of the discretized normal and tangential stresses and the choice of the element (linear Q1 or quadratic Q2) are carried out. The numerical examples illustrate the accuracy and robustness of the proposed mortar-finite element formulation for a contact problem with friction in large deformation.

## 2. Continuous problem and functional framework

One considers the deformation of two elastic bodies occupying, in the initial configuration, two domains $\Omega^j$, j=1,2. For j=1,2 the boundary $\Gamma^j$ of each solid is the union of three non-overlapping parts $\Gamma^j = \Gamma_u^j \cup \Gamma_g^j \cup \Gamma_c^j$. The displacement field is known on $\Gamma_u^j$ (one can suppose, for example, that the $\Omega^j$ solid is embedded in $\Gamma_u^j$). The $\Gamma_g^j$ boundary is submitted to a density of forces noted $g^j \in (L^2(\Omega^j))^2$. Initially the two solids are in contact on the common part of their boundary $\Gamma_c = \Gamma_c^1 = \Gamma_c^2$. The $\Omega^j$ body is submitted to $f^j \in (L^2(\Omega^j))^2$ forces. The normal unit outward vector on $\Omega^j$ is noted $n^j$ as one designated by $\mu \geq 0$ the friction coefficient (supposed constant on $\Gamma_c$ by simplification).

The Coulomb problem of contact with friction consists in finding the $u^j$ displacements and the $\sigma(u^j)$ stresses which verify the following equations and conditions

$$\sigma(u^j) = D^j \varepsilon(u^j) \text{ in } \Omega^j, \qquad [1]$$

$$\text{div}\,\sigma(u^j) + f^j = 0 \text{ in } \Omega^j, \qquad [2]$$

$$\sigma(u^j)n^j = g^j \text{ on } \Gamma_g^j, \qquad [3]$$

$$u^j = 0 \text{ on } \Gamma_u^j, \qquad [4]$$

in which $\varepsilon(u^j)$ represents the linearized strain tensor, $D^j$ is the fourth order tensor satisfying the usual symmetry and ellipticity conditions in elasticity. Equations [1], [2], [3] and [4] respectively designate the behaviour relation, the equilibrium equation and the Neumann and Dirichlet condition.

To introduce the equations onto the $\Gamma_c$ contact zone, the following notations are adopted

$$u^j = (u^j.n^j)n^j + u_t^j = u_n^j + u_t^j, \qquad \sigma(u^j)n^j = \sigma_n(u^j)n^j + \sigma_t(u^j)t, \qquad [5]$$

where $u_n^j$ and $u_t^j$ respectively represent the normal and tangential displacements and $\sigma_n(u^j)$ and $\sigma_t(u^j)$ respectively designate the normal and tangential stresses where t is a unitary fixed tangent vector.

The equations modelling the unilateral contact on $\Gamma_c$ become

$$[u_n] \leq 0, \qquad \sigma_n(u) \leq 0, \qquad \sigma_n(u)[u_n] = 0. \qquad [6]$$

The $[u_n]$ notation represents the $(u^1.n^1 + u^2.n^2)$ jump of the normal displacement through the $\Gamma_c$ contact zone. The [6] conditions, expressing the unilateral contact between the two bodies, describe respectively the non-penetration condition, the sign condition on the normal stress and the complementary condition.

The conditions of Coulomb's friction on $\Gamma_c$ are written

$$\begin{cases} |\sigma_t(u)| \leq \mu|\sigma_n(u)|, \\ |\sigma_t(u)| < \mu|\sigma_n(u)| \Rightarrow [\dot{u}_t] = 0, \\ |\sigma_t(u)| = \mu|\sigma_n(u)| \Rightarrow \exists \lambda \geq 0 \text{ s.t. } [\dot{u}_t] = -\lambda \sigma_t(u). \end{cases} \qquad [7]$$

Here $[\dot{u}_t]$ represents the jump of the tangential velocity through $\Gamma_c$.

Remark: in this paper, we are interested in the discrete formulation with Lagrange multipliers of the friction. Then we restrict ourself to the displacement formulation of the friction and we replace $[\dot{u}_t]$ in [7] by $[u_t]$.

Let us consider K as the closed convex cone of admissible displacements which satisfies the conditions of non penetration

$$K = \{v = (v_1, v_2) \in V^1 \times V^2; [v_n] \leq 0 \text{ on } \Gamma_c\} \qquad [8]$$
$$\text{where } V^j = \{v^j \in (H^1(\Omega^j))^2, v = 0 \text{ on } \Gamma_u^j\}.$$

The variational formulation corresponding to problem [1]-[7], obtained by Duvaut and Lions [20] consists of finding u which verifies

$$u \in K, \ a(u, v-u) + j(u,v) - j(u,u) \geq L(v-u), \ \forall v \in K \qquad [9]$$

where

$$a(u,v) = a^1(u,v) + a^2(u,v), \quad a^j(u,v) = \int_{\Omega^j} (D^j \varepsilon(u^j) : \varepsilon(v^j)) d\Omega, \qquad [10]$$

$$L(v) = \sum_{j=1}^{2} \int_{\Omega^j} f^j \cdot v^j d\Omega + \int_{\Gamma_u^j} g^j \cdot v^j d\Gamma, \qquad j(u,v) = \int_{\Gamma_c} \mu|\sigma_n(u)|[v_t]|d\Gamma, \qquad [11]$$

are defined for all u and v in Sobolev's standard space $(H^1(\Omega))^2$. The functional $j(.,.)$ translates friction.

## 3. Mixed variational formulation of the discrete problem

Let $\mathfrak{I}_h^j$ be a regular family of partitions of $\Omega^j$ into triangles (or quadrangles) $\kappa$

$$\Omega^j = \bigcup_{\kappa \in \mathfrak{S}_h^j} \kappa. \qquad [12]$$

The discretisation parameter $h_j$ on $\Omega^j$ is given by

$$h_j = \max_{\kappa \in \mathfrak{S}_h^j} h_\kappa \qquad [13]$$

where $h_\kappa$ denotes the diameter of the triangle (quadrangle) $\kappa$. Let $h = \max(h_1, h_2)$. For any integer $q \geq 0$, the notation $P_q(\kappa)$ denotes the space of the polynomials with the global degree $\leq q$ on $\kappa$. The finite element space used in $\Omega^j$ is then defined by

$$V_h^j = \left\{ v_h^j \in (C(\overline{\Omega}^j))^2, \forall \kappa \in \mathfrak{S}_h^j, v_h^j \big|_\kappa \in (P_q(\kappa))^2, v_h^j \big|_{\Gamma_u^j} = 0 \right\}, q = 1 \text{ or } 2, \qquad [14]$$

and the approximation space of V becames $V_h = V_h^1 \times V_h^2$.

The contact zone $\Gamma_c$ inherits two independent regular families of monodimensional meshes. The set of nodes belonging to triangulation $\mathfrak{S}_h^j$ are denoted

$$\xi_h^j = \left\{ x_0^j < x_1^j < \ldots < x_{N(h)-1}^j < x_{N(h)}^j \right\}. \qquad [15]$$

In order to express the constraints by using conveniently chosen Lagrange multipliers on the contact zone, we have to introduce first the space describing the degree of the polynomial approximation

$$W_h^{1,j}(\Gamma_c) = \left\{ v_h^j|_{\Gamma_c}, v_h^j \in V_h^j \right\} = W_{hn}^{1,j}(\Gamma_c) \times W_{ht}^{1,j}(\Gamma_c)$$

$$W_h^{0,j}(\Gamma_c) = \left\{ \mu_h, \mu_h \big|_{]z_k^j, z_{k+1}^j[} \in P_0\left(\right]z_k^j, z_{k+1}^j[\right), 0 \leq i \leq N(h) \right\} = W_{hn}^{0,j}(\Gamma_c) \times W_{ht}^{0,j}(\Gamma_c) \qquad [16]$$

where $z_0^j = x_0^j$, $z_{N(h)+1}^j = x_{N(h)}^j$ and for k=1,...,N(h)-1, $z_k^j$ denotes the middle of segment $T_{k-1}^j = \left[x_{k-1}^j, x_K^j\right]$.

The Lagrange multipliers associated to the normal and tangential stresses on the $\Gamma_c$ contact surface either belong to the space $W_h^1(\Gamma_c)$ consisting of continuous piecewise linear functions either belong to the $W_h^0(\Gamma_c)$ space consisting of constant piecewise functions (figure 1).

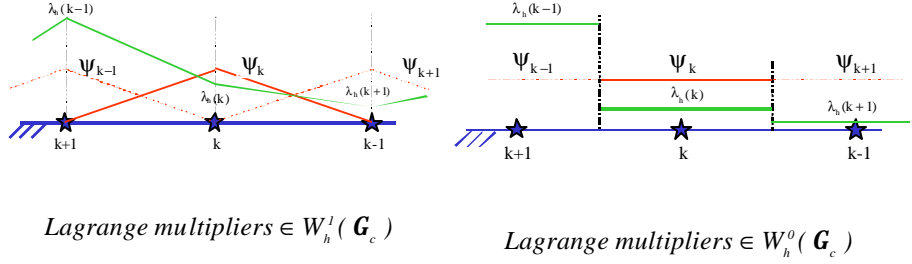

*Lagrange multipliers $\in W_h^1(\mathbf{G}_c)$*     *Lagrange multipliers $\in W_h^0(\mathbf{G}_c)$*

**Figure 1.** *Graphic representation of the two Lagrange multiplier spaces*

Next, we introduce the convex cones associated to the normal and tangential stresses on the contact zone $\Gamma_c$. Let $M_h^1 = M_{hn}^1 \times M_{ht}^1$ be the convex sets of continuous piecewise linear Lagrange multipliers

$$M_{hn}^1 = \left\{ \mu_h \in W_{hn}^1(\Gamma_c), \int_{\Gamma_c} \mu_h \psi_h d\Gamma \geq 0, \ \forall \psi_h \in W_{hn}^1(\Gamma_c), \psi_h \geq 0 \text{ on } \Gamma_c \right\},$$

$$M_{ht}^1 = \left\{ \mu_h \in W_{ht}^1(\Gamma_c), \left| \int_{\Gamma_c} \mu_h \psi_h d\Gamma \right| \leq \int_{\Gamma_c} s_h \psi_h d\Gamma, \ \forall \psi_h \in W_{ht}^1(\Gamma_c) \right\},$$

[17]

where $s_h$ is the given slip bound on $\Gamma_c$. We then consider the convex sets of piecewise constant Lagrange multipliers denoted $M_h^0 = M_{hn}^0 \times M_{ht}^0$ and defined on $\Gamma_c$ as follows

$$M_{hn}^0 = \left\{ \mu_h \in W_{hn}^0(\Gamma_c), \int_{\Gamma_c} \mu_h \psi_h d\Gamma \geq 0, \ \forall \psi_h \in W_{hn}^1(\Gamma_c), \psi_h \geq 0 \text{ on } \Gamma_c \right\},$$

$$M_{ht}^0 = \left\{ \mu_h \in W_{ht}^0(\Gamma_c), \left| \int_{\Gamma_c} \mu_h \psi_h d\Gamma \right| \leq \int_{\Gamma_c} s_h \psi_h d\Gamma, \ \forall \psi_h \in W_{ht}^1(\Gamma_c) \right\}.$$

[18]

In order to solve the Coulomb's frictional contact problem [9] with Lagrange multipliers method, we introduce the following intermediary problem with a given slip limit $s_h$ (see[2,19] for detailed study)

$$P(s_h) \begin{cases} \text{Find } (u_h, \lambda_{hn}, \lambda_{ht}) \in V_h \times M_{hn} \times M_{ht} \text{ such that:} \\ a(u_h, v_h) - \int_{\Gamma_c} \lambda_{hn}[v_{hn}]d\Gamma - \int_{\Gamma_c} \lambda_{ht}[v_{ht}]d\Gamma = L(v_h), \quad \forall v_h \in V_h, \\ \int_{\Gamma_c} (v_{hn} - \lambda_{hn})[u_{hn}]d\Gamma + \int_{\Gamma_c} (v_{ht} - \lambda_{ht})[u_{ht}]d\Gamma \leq 0, \\ \forall (v_{hn}, v_{ht}) \in M_{hn} \times M_{ht}, \end{cases} \quad [19]$$

where $M_{hn} = M_{hn}^1$ or $M_{hn}^0$ and $M_{ht} = M_{ht}^1$ or $M_{ht}^0$.

The discrete mixed problem $P(s_h)$ admits a unique solution (see[21]). It becomes then possible to define a map $\Phi_h$ as follows

$$\Phi_h : M_{hn} \to M_{hn} \qquad [20]$$
$$s_h \to \lambda_{hn}$$

where $(u_h, \lambda_{hn}, \lambda_{ht})$ is the solution of $P(s_h)$. The introduction of this map allows the definition of a discrete solution of Coulomb's frictional contact problem [9].

## 4. Matrix formulation of the global type approach

The matrix formulation of the mixed problem of two bodies $\Omega^1$ and $\Omega^2$ in contact is given by fixing h, the element lengths. One then has a discretization including $N=N_1+N_2$ nodes where $N_1$ is the number of nodes belonging to $\Omega^1$. The N basic functions of $V_h$ are noted $\varphi_i$, $i=1,\ldots,N$ so that if $u_h = (u_h^1, u_h^2)$ we have

$$u^1 = \sum_{i=1}^{N_1} u_h^1(i)\varphi_i \quad \text{and} \quad u^2 = \sum_{i=N_1+1}^{N} u_h^2(i)\varphi_i \, . \qquad [21]$$

We designate by m the number of nodes (i=1,…m) on $\Gamma_c^1$ (slave surface) belonging to the $\Omega^1$ mesh and by n (i= $N_1+1,\ldots N_1+n+1$) the number of nodes on $\Gamma_c^2$ (master surface) belonging to the $\Omega^2$ mesh. The discrete multipliers of the normal and tangential contact stresses are defined on $\Gamma_c^1$ as follows

$$\lambda_{hn} = \sum_{k=1}^{m} \lambda_{hn}(k)\psi_k \quad \text{and} \quad \lambda_{ht} = \sum_{k=1}^{m} \lambda_{ht}(k)\psi_k \, , \qquad [22]$$

where $\psi_k$ are the m basic functions on $\Gamma_c^1$ at the k nodes.

The first discrete formulation equation of the contact problem with friction on the $\Omega^1$ domain has the following matricial form

$$K^1 U^1 - \begin{pmatrix} G_N^1 \\ G_T^1 \\ 0 \end{pmatrix} \Lambda = F^1, \qquad [23]$$

where $K^1$ designates the elastic rigidity matrix linked to $\Omega^1$, $U^1$ designates the vector whose components are the nodal values of $u_h^1$ and $\Lambda = (\Lambda_N, \Lambda_T)$ the vector of components $\lambda_{hn}(k)$, $\lambda_{ht}(k)$ for k=1,…,m. The vector of exterior forces is noted $F^1$ whereas $G_N^1, G_T^1$ are the coupling symmetrical matrices (of order m) between multipliers and displacements. The coefficients of $G_N^1$ and $G_T^1$ matrices are respectively defined by

$$a_{i,j}^1 = \int_{\Gamma_c} \varphi_i . n^1 \psi_j d\Gamma, 1 \leq i, j \leq m,$$
$$a_{i,j}^1 = \int_{\Gamma_c} \varphi_i . t^1 \psi_j d\Gamma, 1 \leq i, j \leq m. \qquad [24]$$

Remark : for a fixed choice of all the multipliers, the $G_N^1$ and $G_T^1$ matrices are identical and will be noted $G^1 = G_N^1 = G_T^1$. In the same manner, the system of unknown equations $U^2$ and $\Lambda$ on $\Omega^2$ is written

$$K^2 U^2 - \begin{pmatrix} G_N^{2,1} \\ G_T^{2,1} \\ 0 \end{pmatrix} \Lambda = F^2, \qquad [25]$$

$G^{2,1} = G_N^{2,1} = G_T^{2,1}$ is a rectangular matrix of n lines and m columns whose coefficients are

$$a_{i,j}^{2,1} = \int_{\Gamma_c} \varphi_i . n^2 \psi_j d\Gamma, \; N_1 \leq i \leq N_1 + n + 1 \text{ et } 1 \leq j \leq m. \qquad [26]$$

Finally the problem of contact with friction between the two bodies is written

$$\begin{pmatrix} K^1 & 0 \\ 0 & K^2 \end{pmatrix} \begin{pmatrix} U^1 \\ U^2 \end{pmatrix} - \begin{pmatrix} G^1 \\ 0 \\ G^{2,1} \\ 0 \end{pmatrix} \Lambda = \begin{pmatrix} F^1 \\ F^2 \end{pmatrix}. \qquad [27]$$

The interest now is in the matricial writing of the contact and friction conditions

$$\int_{\Gamma_c} (v_{hn} - \lambda_{hn})[u_{hn}]d\Gamma + \int_{\Gamma_c} (v_{ht} - \lambda_{ht})[u_{ht}]d\Gamma \le 0, \quad \forall (v_{hn}, v_{ht}) \in M_h, \qquad [28]$$

with $M_h = M_{hn}^1 \times M_{ht}^1$ or $M_h = M_{hn}^0 \times M_{ht}^0$. Let $U_N$ and $U_T$ be the vectors whose components are respectively the nodal values of $[u_{hn}]$ and $[u_{ht}]$. It can be shown [19] that the preceding inequation is written

$$\begin{cases} (G^1 \Lambda_N)_i & \le 0, \\ (U_N^1 + (G^1)^{-1} - (G^{2,1})^t U_N^2)_i & \le 0, \\ (G^1 \Lambda_N)_i (U_N^1 + (G^1)^{-1} - (G^{2,1})^t U_N^2)_i \le 0, \\ \left| (G^1 \Lambda_T)_i \right| & \le -\mu(\Lambda_N)_i, \\ \left| (G^1 \Lambda_T)_i \right| < -\mu(\Lambda_N)_i & \Rightarrow (U_T^1 + (G^1)^{-1} - (G^{2,1})^t U_T^2)_i = 0, \\ (G^1 \Lambda_T)_i (U_T^1 + (G^1)^{-1} - (G^{2,1})^t U_T^2)_i & \le 0. \end{cases} \qquad [29]$$

### 4.1 Construction of the $G^1$ and $G^{2,1}$ matrices

In the case of the basic functions $\psi_k$ on $\Gamma_c^1$, continuous piecewise linear (P1) or constant piecewise (P0) and for the Q1 finite elements, the construction of the $G^1$ and $G^{2,1}$ matrices is described in this paragraph. The $G^1$ matrix coefficients for functions $\psi_k$ of the P0 and P1 type are shown on figure 2 and 3 respectively.

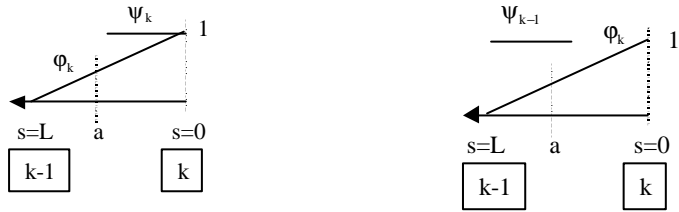

$$a_{k,k} = \int_0^a \varphi_k \psi_k ds = a(-\frac{a}{2L}+1) \quad, a \leq \frac{L}{2}$$
$$a_{k,k} = \frac{3}{8}L \quad, a > \frac{L}{2}$$

$$a_{k,k-1} = \int_a^{\frac{L}{2}} \varphi_k \psi_{k-1} ds = \frac{a^2}{2L} \quad, a \geq \frac{L}{2}$$
$$a_{k,k-1} = \frac{L}{8} \quad, a < \frac{L}{2}$$

**Figure 2.** *Coefficient of $G^1$ for P0 shape functions $\psi_k$*

The $G^{2,1}$ matrix is the matrix coupling the m slave nodes of the $\Gamma_c^1$ surface and the n master nodes of the $\Gamma_c^2$ surface. To determine the expression of the coefficients of this coupling matrix, in the case when the surfaces are not smooth (figure 4a), it is necessary to proceed first of all to a projection of the interface nodes onto a curvilinear abscissa that is noted "s" (figure 4.b).

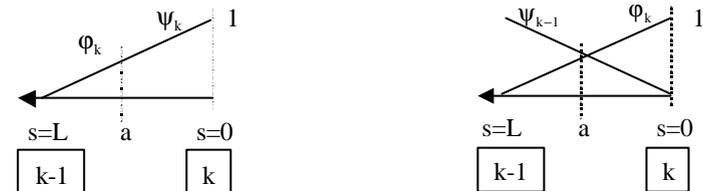

$$a_{k,k} = \int_0^a \varphi_k \psi_k ds = a(\frac{a^2}{3L^2} - \frac{a}{L} + 1)$$
$$a_{k,k-1} = \int_0^a \varphi_k \psi_{k-1} ds = a(-\frac{a^2}{3L^2} + \frac{a}{2L})$$

**Figure 3.** *Coefficient of $G^1$ for P1 shape functions $\psi_k$*

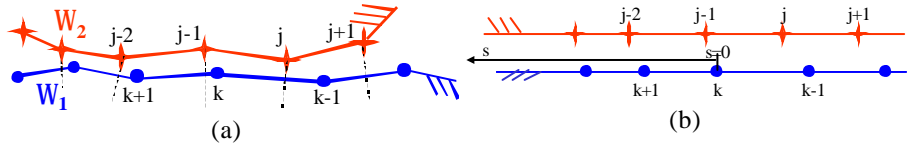

**Figure 4.** *a. Contact surfaces $G_c^1$ and $G_c^2$ at time t ; b. Projection of the contact surfaces $G_c^1$ and $G_c^2$ on the curvilinear abscissa s*

If one wishes to fill the $G^{2,1}$ matrix column corresponding to the k slave node, one calculates the curvilinear abscissa of the nodes of $\Gamma_c^1$ and $\Gamma_c^2$ by fixing the origin s=0 to the node k.

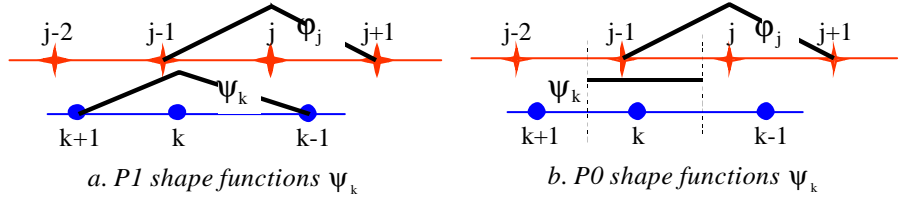

*a. P1 shape functions $\psi_k$*     *b. P0 shape functions $\psi_k$*

**Figure 5.** *Coefficient calculation of the matrix $G^{2,1}$ for the slave node k*

The non null coefficients of $G^{2,1}$ for the k slave node are $a^{2,1}_{j-2,k}, a^{2,1}_{j-1,k}, a^{2,1}_{j,k}, a^{2,1}_{j+1,k}$ or $a^{2,1}_{j-2,k}, a^{2,1}_{j-1,k}, a^{2,1}_{jk}$ if one chooses the $\psi_k$ basis functions of the P1 type (figure 5.a) or of the P0 type (figure 5.b).

## 5. Numerical Results

In this section, one studies and compares numerically the performances of the methods shown previously in the case of contact with friction or without friction ; the analysis of the quality of approximation of these methods having been presented in[22,2,19]. These methods have been implemented into PLAST2[1], a finite elements code in explicit dynamics based on the method of the Lagrange multipliers. This code deals with contact and friction conditions with either the Lagrange interpolation operator (local type approach) or the mortar-finite element approach (global type approach). For the first approach, contact is defined for each node of the slave surface by using the intervention of the closest segment defined by 2 nodes of the master surface. This gives to the condition (also called node-on-segment contact condition), a very local characteristic observed on the different chosen tests.

### 5.1. First numerical test

In this numerical test, one considers the contact problem shown in Figure 6. The $\Omega^1$ domain is a part of a disc of 1 mm radius, the $\Omega^1$ domain is a rectangle of 1.8 mm x 0.3 mm. On each domain, the behavior law is that of Hooke for the isotropic and homogeneous materials. For k=1,2

$$\sigma^k_{ij}(u^k) = \frac{E_k \nu_k}{(1-2\nu_k)(1+\nu_k)} \delta_{ij} \varepsilon^k_{mm}(u^k) + \frac{E_k}{1+\nu_k} \varepsilon^k_{ij}(u^k), \qquad [30]$$

with $E_1$=70000MPa, $E_2$=7000MPa and $\nu_1=\nu_2$=0.3.

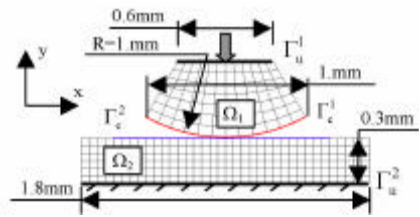
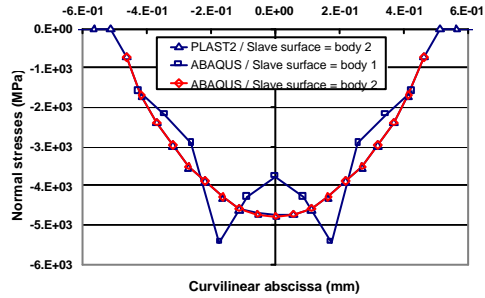

**Figure 6.** *Finite element model*      **Figure 7.** *Normal stresses with the local type approach*

The $\Omega^1$ domain is embedded on the $\Gamma_u^2$ boundary. The displacement cycle imposed on $\Gamma_u^1$ is a vertical indentation of -0.1581mm.

On each solid, one uses rectangular finite elements of the Q1 type (4 node quadrangles) or of the Q2 quadratic type (8 node quadrangles) [23]. Let us note that the normal and tangential stresses are represented on each figure for a maximum indentation of -0.1581mm.

*5.1.1. The local type approach: comparison of PLAST2 and ABAQUS_Standard codes*

The aim is to compare PLAST2 and ABAQUS codes on a problem of contact without friction using the classical node-on-segment approach. This allows us firstly to validate the PLAST2 code and to show the limits of the local type approach for dealing with the contact conditions. Let us note that to solve problem [27] using a dynamic code, one replaces the displacement cycle imposed on $\Gamma_u^1$ by a very weak vertical speed (damping and inertia terms are therefore negligible) subjected to the same surface that puts the two solids under the same deformation cycle. Since the contact surface deforms during this cycle, it is necessary to update the $G^1$ and $G^{2,1}$ coupling matrixes at each time step increment.

Figure 7 represents the distribution of the contact normal stresses for a maximum indentation when $\Gamma_c^1$ or $\Gamma_c^2$ is the slave surface. One will note the similarities of the stresses calculated by the two codes and the asymmetrical results obtained when the slave surface is changed. It is clear that on this test the local type approach has shown its limits.

*5.1.2 The global type approach in PLAST2 with Q1 finite elements type*

In this paragraph, another technique to approximate contact problems implemented in PLAST2 is presented, called global type approach. One considers first of all the problem of contact without friction. This involves studying the behavior of the global type approach whether $\Gamma_c^1$ or $\Gamma_c^2$ is chosen as a slave surface.

By considering the case of piecewise constant multipliers on the contact interface ($M_h^0$), the symmetrical behavior of the global type approach compared to the local one is observed on the distribution of the contact stresses (Figure 8) obtained for the maximum indentation and for the choice of whatever slave surface ($\Gamma_c^1$ or $\Gamma_c^2$).

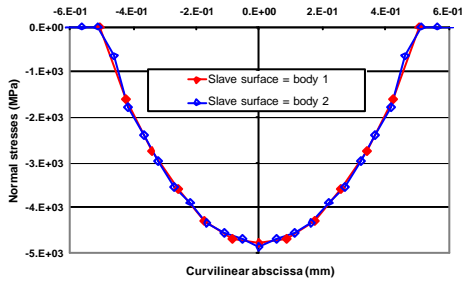
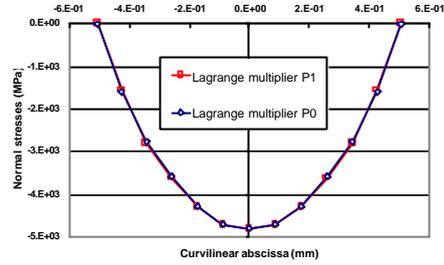

**Figure 8.** *Normal contact stresses (maximum indentation) when $G_c^1$ or $G_c^2$ is the slave surface with the global type approach*

**Figure 9.** *Normal contact stresses ($M_{hn} = M_{hn}^0$ or $M_{hn} = M_{hn}^1$) for $G_c^1$ slave surface*

The comparison of figures 7 and 8 shows that the global type approach makes the management of the contact more symmetrical when slave surfaces are interchanged.

On figure 9, one will note the similarities of the contact normal stresses when $M_{hn} = M_{hn}^0$ and $M_{hn} = M_{hn}^1$ (P0 and P1 multipliers respectively).

Let us now consider the case of a problem of contact with Coulomb's friction. One has to insure the correct behavior of the global type approach by using the different convex approximations ($M_h = M_h^0$ or $M_h = M_h^1$) linked to normal and tangential stresses. In this case, the preceding comments apply. In particular, the shape of the normal and tangential stresses is similar (Figure 10) for the two types of approximations of $M_h$ and for different µ friction coefficients.

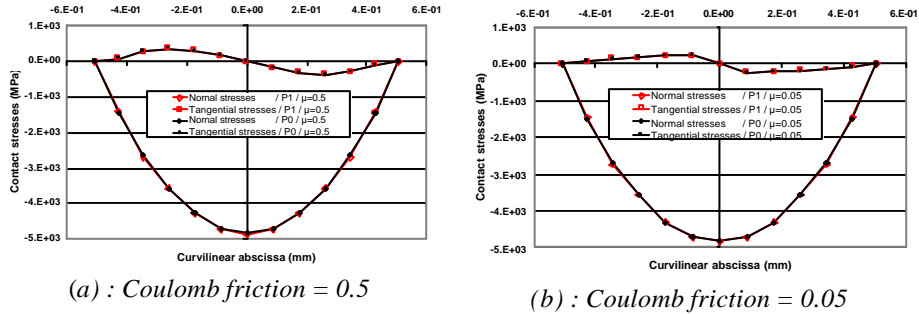

(*a*) : *Coulomb friction = 0.5*  (*b*) : *Coulomb friction = 0.05*

**Figure 10.** *Normal and tangential stresses ( $M_{hn} = M_{hn}^0$ or $M_{hn} = M_{hn}^1$ ) for different Coulomb friction*

### 5.1.3. The global type approach in PLAST2 with Q2 finite elements type

The global type approach has been implemented into PLAST2 for quadratic finite elements Q2 to simulate the problem of contact with Coulomb's friction between two elastic solids (see [18] for a problem of unilateral contact without friction).

The convexes of Lagrange multipliers are continuous piecewise linear functions ($M_h^1$) or constant piecewise functions ($M_h^0$) on $\Gamma_c$. The use of such finite elements gives hope for a better precision of calculation compared to rectangular or linear elements[24]. Figures 11 and 12 validate the correct behavior of the global type approach for problems of contact with or without friction. In the same manner for finite elements of the Q1 type, the results are identical for the Lagrange multipliers $M_h = M_h^0$ or $M_h = M_h^1$.

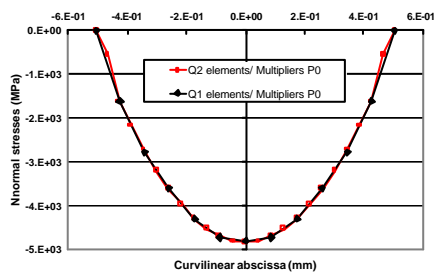
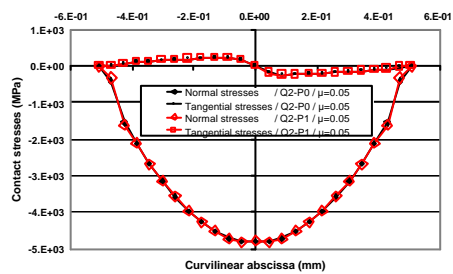

**Figure 11.** *Normal stresses ( $M_{hn} = M_{hn}^0$ ) for Q1 and Q2 element*

**Figure 12.** *Normal and tangential stresses ( $M_{hn} = M_{hn}^0$ or $M_{hn} = M_{hn}^1$ ) for Q2 element type and Coulomb*

*type without friction*          *friction m=0.05*

### 5.2. Second numerical test

In the case of contact with friction of a deformable body on a rigid surface, one studies numerically the performance of the methods shown above for $M_h = M_h^0$ and $M_h = M_h^1$. The numerical tests have been carried out on the finite element code PLAST2. In the numerical tests, the behavior law is that of Hooke for isotropic and homogeneous materials

$$\sigma_{ij}(u) = \frac{E\nu}{(1-2\nu)(1+\nu)} \delta_{ij} \varepsilon_{mm}(u) + \frac{E}{1+\nu} \varepsilon_{ij}(u), \qquad [31]$$

with $E=7.10^4$ MPa et $\nu=0.3$.

The $\Omega$ domain is a rectangle measuring 1.3 mm x 0.3 mm. The discretisation is carried out with finite rectangular elements of the Q1 type in plane strains. The origin of the curvilinear abscissa is defined from the point O in the trigonometric direction. A total displacement of $2.10^{-3}$ mm is imposed on the $\Gamma_u^1$ and $\Gamma_u^2$ (see Figure 13). The horizontal displacement is null on $\Gamma_u^1$. The vertical displacement on $\Gamma_u^2$ is free which enables a detachment of the deformable body for a curvilinear abscissa superior to 0.7mm (see Figure 15.a). The Tresca threshold stress is equal to 200MPa.

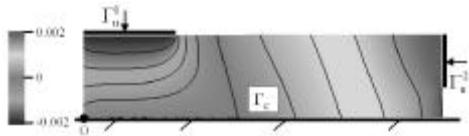
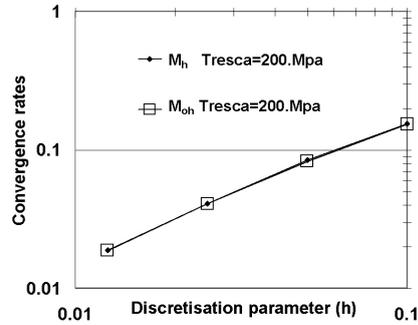

**Figure 13.** *Vertical displacement on the reference mesh*

**Figure 14**. *Convergence rates of the two approach $M_{0h} = M_h^0$ and $M_h = M_h^1$*

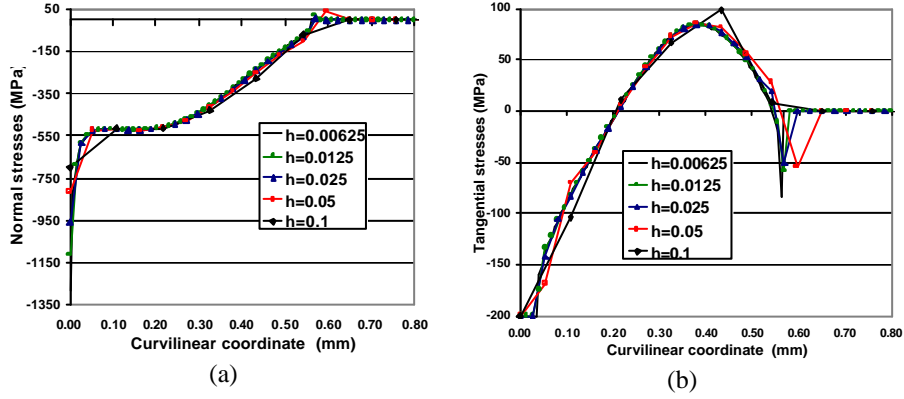

**Figure 15.** *a. Normal ; b. tangential stresses for various h and for* $M_h = M_h^0$

Having no analytic solution for the problem treated, the $\|u - u_h\|$ error in the energy is numerically estimated by $\|u_{ref} - u_h\|$. The reference solution is calculated on a reference mesh containing 9678 elements.

On Figure 14 the convergence order of the different methods for different discretisation parameter h is represented. It can be seen that the convergence is similar for the two approaches. On Figure 15, one can see that the normal stress is not a negative function over all the interface, this is due to the use of slightly negative Lagrange multipliers. However this method enables the singularities of the stresses edges to be attenuated.

### 5.4. Third numerical test

This test enables the global type approach to be validated when compared to the local one, in the case of a simulation of three deformable bodies in contact with one another and with a rigid surface (Figure 16). The bodies have an elasto-plastic behaviour ($\sigma_{eq} = 348(.2\,10^{-2} + \varepsilon^p)^{0.03}$). Penetrations of the master nodes into the slave surface appear in the simulation using the local type approach and they generate a divergence of the calculation whereas the simulation with the global one is carried out without problems.

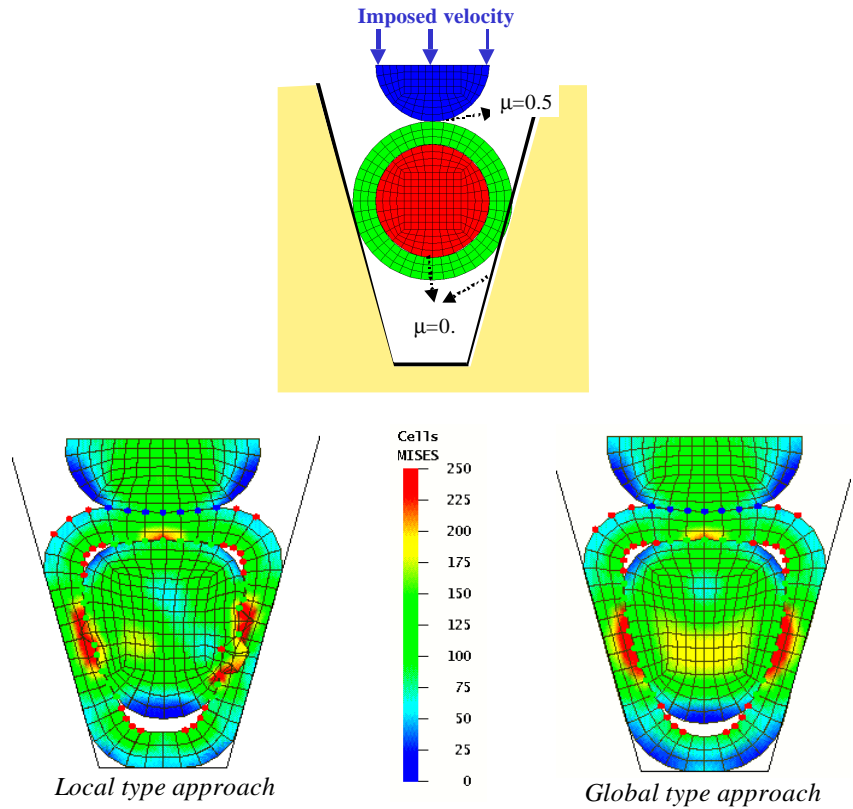

**Figure 16.** *Comparison between local and global type approach on a forging simulation of three deformable bodies*

## 6. Conclusion

For the management of a problem of friction contact with the mortar-finite element method and for the construction of a matrix expressing the friction contact of deformable bodies, one can choose form functions (linked to Lagrange multipliers) which are continuous piecewise linear (P1) or constant piecewise (P0). These two approaches have been implemented in PLAST2. Mathematically and numerically on problems of contact with friction, the results of normal and tangential contact stresses are similar when using both approximations (P0 or P1) of Lagrange multipliers.

For a discretisation with quadrilateral elements of the Q1 or Q2 type, it has been established on various numerical tests, that the mortar-finite element method makes the management of the contact much more symmetrical when slave surfaces are exchanged, thus closer to physics. Finally, it has been shown that the formulation of this global type approach is also appropriate for Q1 or Q2 finite elements.